\documentclass[11pt,a4paper]{amsart}

\usepackage{amsfonts}
\usepackage{amsmath}
\usepackage{amssymb}
\usepackage{amsthm}
\usepackage{enumitem}
\usepackage{geometry}
\usepackage{latexsym}
\usepackage{hyperref}
\usepackage{mathrsfs}
\usepackage{tikz}
\usepackage{tikz-cd}
\usepackage{xy}
\usepackage{color}
\usepackage{verbatim}
\allowdisplaybreaks[4]

\hypersetup{colorlinks=true,citecolor=red,linkcolor=blue,urlcolor=black,}

\usepackage[left,pagewise,mathlines]{lineno}
\newtheorem{theorem}{Theorem}[section]
\newtheorem{lemma}[theorem]{Lemma}
\newtheorem{corollary}[theorem]{Corollary}
\newtheorem{proposition}[theorem]{Proposition}

\newtheorem{def-thm}[theorem]{Definition-Theorem}

\theoremstyle{definition}
\newtheorem{definition}[theorem]{Definition}

\newtheorem{remark}[theorem]{Remark}

\DeclareMathOperator{\codim}{codim}

\DeclareMathOperator{\Spec}{Spec}

\DeclareMathOperator{\Supp}{Supp}

\DeclareMathOperator{\Exc}{Exc}

\DeclareMathOperator{\id}{id}

\DeclareMathOperator{\ord}{ord}
\DeclareMathOperator{\Sing}{Sing}

\DeclareMathOperator{\NE}{\overline{NE}} 
\DeclareMathOperator{\Nlc}{Nlc} 
\DeclareMathOperator{\rank}{rank}
\DeclareMathOperator{\Null}{Null}
\DeclareMathOperator{\ninv}{ninv}
\DeclareMathOperator{\inv}{inv}
\DeclareMathOperator{\lcm}{lcm}
\DeclareMathOperator{\Bs}{Bs}
\DeclareMathOperator{\Der}{Der}
\DeclareMathOperator{\Bl}{Bl}

\newcommand{\cO}{\ensuremath{\mathcal{O}}}
\newcommand{\cF}{\ensuremath{\mathcal{F}}}
\newcommand{\cG}{\ensuremath{\mathcal{G}}}
\newcommand{\Q}{\ensuremath{\mathbb{Q}}}
\newcommand{\R}{\ensuremath{\mathbb{R}}}
\newcommand{\C}{\ensuremath{\mathbb{C}}}

\newcommand{\bN}{\ensuremath{\mathbf{N}}}
\newcommand{\bM}{\ensuremath{\mathbf{M}}}
\newcommand{\bB}{\ensuremath{\mathbf{B}}}
\newcommand{\bir}{\ensuremath{\dashrightarrow}}

\newcommand{\fm}{\ensuremath{\mathfrak{m}}}

\title{On the generalized foliated cone theorem of rank one}
\author{Mengchu Li}

\address{School of Mathematical Sciences, Zhejiang University, Hangzhou 310058, China}
\email{mengchuli@zju.edu.cn}

\subjclass[2020]{14E30, 37F75}
\keywords{Foliations, Generalized pairs, Cone theorem.}

\begin{document}

\begin{abstract}
We prove the cone theorem for generalized foliated quadruples of rank one in arbitrary dimension. As an application, we prove the uniqueness of MMP outputs for log canonical generalized foliated quadruples of rank one.
\end{abstract}

\maketitle

\setcounter{tocdepth}{1}
\tableofcontents

\section{Introduction}

The cone theorem is one of the basic parts of the minimal model program (MMP). It describes negative extremal rays of the cone of curves in terms of rational curves. The minimal model program can also be studied for foliations, with the foliated cone theorem as a basic ingredient. For technical reasons, one also considers generalized foliated quadruples $(X,\cF,B,\bM)$ and studies  their minimal model program. 

For foliations of rank one, Bogomolov and McQuillan \cite[Theorem 4.2.1]{BM16} proved the foliated cone theorem for Gorenstein canonical foliations of rank one. McQuillan \cite{McQ24} proved a stack-theoretic cone theorem, also treating orbifold and logarithmic boundaries. Using foliation adjunction, Cascini and Spicer \cite[Theorem 4.8]{CS25a} proved a cone theorem for rank one foliated triples. The cone theorem for generalized foliated threefolds of rank one was proved in \cite[Theorem 1.2]{Li25}. The foliated cone theorem and minimal model program are also well-developed in various other settings, see \cite{Spi20,CS21,ACSS22,CHLX23,cascini_minimal_2024,CM24,fujino_toric_2025,CZ26,CS26} and the references therein.

Our main result is the following generalized foliated cone theorem of rank one in arbitrary dimension, which generalizes \cite[Theorem 4.8]{CS25a} and \cite[Theorem 1.2]{Li25}.

\begin{theorem}\label{MainThm:ConeTHM}
    Let $X$ be a normal projective variety, and let $(X,\cF,B,\bM)$ be a generalized foliated quadruple of rank 1 on $X$. Assume that $K_\cF$ is $\Q$-Cartier and $B$ is $\R$-Cartier. Then there exists an at most countable family of rational curves $\{C_j\}_{j\in \Lambda}$ tangent to $\cF$ satisfying
    \[
    0 < -(K_\cF + B + \bM_X) \cdot C_j \leq 2
    \]
    for each $j \in \Lambda$. Moreover, 
    \[
    \NE(X) = \NE(X)_{K_\cF + B + \bM_X \geq 0} + \NE(X)_{-\infty} + \sum_{j \in \Lambda} \R_{\geq 0}[C_j],
    \]
    where $\NE(X)_{-\infty}$ is the subcone of $\NE(X)$ spanned by the images of $\NE(W)\to \NE(X)$ where $W$ is a non-lc center of $(X,\cF,B,\bM)$.
\end{theorem}

We briefly explain the proof of Theorem~\ref{MainThm:ConeTHM}. The proof follows the strategy of \cite[Theorem 4.8]{CS25a}. Consider an lc gfq $(X,\cF,B,\bM)$ of rank one, and set $K=K_\cF+B+\bM_X$. Let $H_R$ be a supporting function of a $K$-negative exposed ray $R$, and let $W$ be an irreducible component of $\Null H_R$ not contained in the singular locus of $\cF$. The key step is to show that $W$ is $\cF$-invariant. If $W$ is not $\cF$-invariant, then after replacing $X$ by a suitable birational model, we may assume there exists a non-invariant $K$-negative curve $C$ contained in $W$ such that $\cF$ is terminal at every closed point in an analytic neighborhood of $C$, and there is an $\cF$-invariant surface $\Gamma$ containing $C$. In \cite[Theorem 4.8]{CS25a}, adjunction to $\Gamma$ and then to $C$ gives a contradiction. 
As noted in \cite[Remark 4.9]{Li25}, this step does not directly extend to generalized foliated structures, since it is not clear how to define restricted b-divisors along higher codimensional subvarieties. 
In our setting, Proposition~\ref{prop:extraction} instead provides a  birational model $h:Y\to X$ with an invariant prime divisor $S$ on $Y$ centered on $\Gamma$. We may apply divisorial adjunction to the crepant pullback of $K_\cF + B + \bM_X$ along $S$. Proposition~\ref{prop:extraction} also ensures that the foliation $\cF_Y = h^{-1}\cF$ on $Y$ and the restricted foliation $\cF_S$ on $S^\nu$ are terminal at every closed point. Choosing a general curve $C_S$ over $C$ and an invariant surface containing $C_S$, we may restrict all the data on $S$ and apply induction on dimension; see the proof of Proposition~\ref{prop:nulllocus}. This avoids restricting the b-divisor directly to $\Gamma$ and overcomes the difficulty described in \cite[Remark 4.9]{Li25}.

Another ingredient in the proof of Theorem~\ref{MainThm:ConeTHM} is the Zariski closedness of the non-lc locus of $(X,\cF,B,\bM)$. To prove this, using the resolution algorithm given by Belotto \cite{Bel16}, we construct a log smooth birational model for an lc Gorenstein foliation of rank one on which the induced foliation is also lc. This is stated in Proposition~\ref{prop:lcmodel}, which is also of independent interest.

As an application of Theorem~\ref{MainThm:ConeTHM}, we prove the uniqueness of MMP outputs, starting from a fixed lc generalized foliated quadruple $(X,\cF,B,\bM)$ of rank one, which generalizes \cite[Theorem 1.1]{CMS26}.

\begin{theorem}\label{MainThm:MMPoutputs}
    Let $X$ be a normal projective $\Q$-factorial variety, and let $(X,\cF,B,\bM)$ be a log canonical generalized foliated quadruple of rank 1 on $X$. Assume that $K_\cF + B + \bM_X$ is pseudo-effective. Let $f_i: X \bir Y_i$, $i= 1,2$ be two outputs of $(K_\cF + B + \bM_X)$-MMPs. Then there exists an isomorphism $\phi: Y_1\to Y_2$ such that $f_2 = \phi \circ f_1$.
\end{theorem}

The proof of \cite[Theorem 1.1]{CMS26} relies on a Bertini-type statement \cite[Theorem 3.3]{CMS26}  for invariant subvarieties, using the theory of foliated jets. Our proof of Theorem~\ref{MainThm:MMPoutputs} takes an alternative approach using generalized foliated structures, which allow us to perturb the nef part by a small ample divisor while preserving the log canonical condition. It is worth noting that the uniqueness holds only for MMP outputs from a fixed start $(X,\cF,B,\bM)$. Even for a fixed $(X,\cF)$, its $K_\cF$-minimal models need not be unique up to isomorphism;  see \cite[Remark 3.4]{CMS26} and Remark~\ref{rmk:MMPoutputs}.

\subsection*{Acknowledgement} The author would like to thank Lingyao Xie for useful comments and discussions.

\subsection*{AI disclosure}
The author used an LLM to assist in developing the proof. It drew the author's attention to Belotto's result \cite{Bel16}, which the author used in the proof of Proposition~\ref{prop:lcmodel}. It also suggested the stack-theoretic proof of Proposition~\ref{prop:extraction} and gave the example in Remark~\ref{rmk:HirSur}. The author also used the LLM to check the arguments and the grammar. The author takes full responsibility for the accuracy and content of the paper.

\section{Preliminaries}
We work over the field of complex numbers $\C$. All varieties are assumed to be quasi-projective. We adopt standard notation and definitions in birational geometry as in \cite{KM98} and \cite{BCHM10}, and use them freely in this paper. Throughout this paper, a \emph{subvariety} means an integral closed subscheme.

\subsection{Divisors} 
Let $M$ be a Weil $\R$-divisor on a normal variety $X$.  Denote by $\mu_DM$ the coefficient of a prime divisor $D$ in $M$. Denote by $M^+, M^-$ the positive part and negative part of $M$, i.e. $M = M^+ - M^-$ where $M^+,M^-$ are effective $\R$-divisors with no common components. 

Let $X$ be a normal projective variety and let $L$ be an $\R$-Cartier divisor on $X$. Define the \emph{exceptional locus} of $L$ to be 
\[
\Null L = \bigcup_{L|_V \text{ is not big}} V,
\]
where the union runs over the subvarieties $V$ of $X$ with positive dimension. If $L$ is nef, then $\Null L$ is called the \emph{null locus} of $L$. Let $M$ be a $\Q$-Cartier divisor on $X$. We define the \emph{stable base locus} of $M$ to be 
\[
\mathbf{B}(M)  = \bigcap_{m \in \mathbb{N},\, mM \text{ Cartier}} \Bs \,\lvert mM \rvert.
\]
The \emph{augmented base locus} of $L$ is defined by 
\[
\mathbf{B}_+(L) = \bigcap \mathbf{B}(L-A),
\]
where $A$ runs over all ample $\R$-Cartier divisors on $X$ such that $L-A$ is $\Q$-Cartier. By \cite[Theorem 1.4]{birkar_augmented_2017}, for any nef $\R$-divisor $H$ on $X$, $\Null H = \mathbf{B}_+(H)$. 

\subsection{Models} Let $X$ be a normal projective variety and let $D$ be an $\R$-Cartier $\R$-divisor on $X$. Let $\phi: X \bir X'$ be a birational map. We say that $\phi$ is \emph{$D$-negative} if the following conditions hold:
\begin{enumerate}
    \item $\phi$ is a birational contraction, i.e. $\phi^{-1}$ does not contract any divisors.
    \item $D'= \phi_* D$ is $\R$-Cartier.
    \item For some resolution of indeterminacy $p: W\to X$ and $q: W\to X'$, we may write
    \[
    p^*D = q^*D' + E,
    \]
    where $E\geq 0 $ is $q$-exceptional and $\Supp E$ contains the strict transform of all $\phi$-exceptional divisors. 
\end{enumerate}
If $\phi: X \bir X'$ is $D$-negative and $D' = \phi_* D$ is nef, then $\phi$ is said to be a \emph{minimal model} of $D$, or a $D$-minimal model.

\subsection{Foliation}

A \emph{foliation} of rank $r$ on a normal variety $X$ is a coherent subsheaf $T_\cF\subset T_X$ of $\rank r$ such that $T_X/T_\cF$ is torsion-free and $T_\cF$ is closed under the Lie bracket. The canonical divisor of $\cF$ is a Weil divisor on $X$ defined by $\cO(-K_\cF) \cong \det (T_\cF)$. We may define the \emph{Pfaff field} associated to $\cF$ to be the morphism
\[
\phi: \Omega_X^{[r]} \to \cO_X(K_\cF)
\]
by taking the double dual of the $r$th exterior power of the map $\Omega_X^{1}\to T_\cF^*$, induced by $T_\cF \subset T_X$. Then we have the induced morphism
\[
\phi': \left(\Omega_X^{[r]} \otimes \cO(-K_\cF)\right)^{**} \to \cO_X,
\] 
and we define the \emph{singular locus} $\Sing \cF$ of $\cF$ to be the cosupport of the image ideal of $\phi'$.

Let $X$ be a normal variety and $\cF$ be a foliation on $X$. Let $f: Y \bir X$ be a dominant map. Denote by $f^{-1}\cF$ the \emph{induced foliation} of $\cF$ on $Y$ as in \cite[3.2]{Dru21}. If there exists a dominant map $g: X \bir Z$ such that $\cF = g^{-1}\cF_Z$, where $\cF_Z$ is the foliation by points on $Z$, then we say $\cF$ is \emph{algebraically integrable}. If $h: X\bir X'$ is a birational map, then we define the pushforward foliation $h_*\cF$ of $\cF$ to $X'$ to be $(h^{-1} )^{-1}\cF$. 

We say a subvariety $S\subset X$ is \emph{$\cF$-invariant} if $T_\cF$ is locally free at the generic point of $S$, and for an affine open neighborhood $U$ of the generic point of $S$, $\partial(I_{S\cap U})\subset I_{S\cap U}$ for any $\partial \in H^0(U,T_\cF|_U)$, where $I_{S\cap U}$ is the ideal sheaf of $S \cap U$. We say a subvariety $C$ is \emph{tangent} to $\cF$ if $T_C \subset T_\cF|_C$ at general points of $C$. For a prime divisor $D$ on $X$, we define $\epsilon (D) = 0$ if $D$ is $\cF$-invariant, and $\epsilon (D) = 1$ otherwise. For an $\R$-divisor $M$ on $X$, we define
\[
M^{\ninv } = \sum_{D} \epsilon(D) \mu_D (M) D
\] 
as the \emph{non-invariant part} of $M$, and define $M^{\inv} = M - M^{\ninv}$ as the \emph{invariant part} of $M$.

We introduce a definition of $\cF$-invariant subvarieties on a birational model over $X$.
\begin{definition}
    Let $X$ be a normal variety and let $\cF$ be a foliation on $X$. Let $f: Y \to X$ be a birational morphism. Let $S$ be a subvariety on $Y$. Let $\eta_S$ be the generic point of $S$. We say that $S$ is \emph{$\cF$-invariant over $X$}, if there exists an affine open neighborhood $U$ of $f(\eta_S)$ such that $T_\cF|_U$ is free and an open neighborhood $V\subset f^{-1}(U)$ of $\eta_S$, such that $\partial(I_{S \cap V}) \subset I_{S\cap V}$ for any $\partial \in H^0(U,T_\cF|_U)$, where $I_{S\cap V}$ is the ideal sheaf of $S\cap V$, and we regard $\partial$ as a derivation of the function field $K(X) = K(Y)$.
\end{definition}

It is clear that for a subvariety $S$ on $X$, $S$ is $\cF$-invariant if and only if $S$ is $\cF$-invariant over $X$. However, if $f:Y \to X$ is a birational morphism between normal varieties and $\cF_Y = f^{-1}\cF$, then
$\cF$-invariance over $X$ and $\cF_Y$-invariance of subvarieties of $Y$ do not coincide in general. We introduce this notion to describe the invariance of the blowing-up centers arising in functorial resolutions of $X$.

\subsection{Generalized foliated structure} A \emph{generalized foliated quadruple} $(X,\cF, B,\bM)$ (\emph{gfq} for short) consists of a normal variety $X$, a foliation $\cF$ on $X$, an $\R$-divisor $B\geq 0$ on $X$, a b-nef b-$\R$-Cartier b-divisor $\bM$ on $X$, such that $K_\cF + B + \bM_X$ is $\R$-Cartier. For the definition of a b-divisor and related terminology, see \cite[Definition~2.4]{Li25}. Let $E$ be a prime divisor over $X$, and let $f: Y \to X$ be a birational morphism such that $E$ is a prime divisor on $Y$. We write
\[K_{\cF_Y} + B_Y + \bM_Y = f^*(K_\cF + B + \bM_X),\]
where $\cF_Y = f^{-1}\cF$, and define the \emph{discrepancy} $a(E,\cF,B,\bM) = -\mu_E B_Y$. We say that $(X,\cF,B,\bM)$ is \emph{log canonical} (\emph{lc} for short) if $a(E,\cF,B,\bM) \geq -\epsilon (E)$ for every prime divisor $E$ over $X$. We say that $(X,\cF,B,\bM)$ is \emph{canonical} (resp. \emph{terminal}) if  $a(E,\cF,B,\bM) \geq 0$ (resp. $>0$) for any exceptional divisor $E$ over $X$. We use the prefix ``sub-'' if $B$ is not required to be effective. An \emph{lc place} (resp. \emph{non-lc place}) of $(X,\cF,B,\bM)$ is a prime divisor $E$  over $X$ such that $a(E,\cF,B,\bM) = -\epsilon (E)$ (resp. $<-\epsilon(E)$). An \emph{lc center} (resp. \emph{non-lc center}) is the center of an lc place (resp. non-lc place) on $X$. The union $\Nlc(X,\cF,B,\bM)$ of non-lc centers of $(X,\cF,B,\bM)$ is called the \emph{non-lc locus} of $(X,\cF,B,\bM)$. 

For any $P\in X$ not necessarily closed, we say $(X,\cF,B,\bM)$ is terminal (resp. canonical) at $P$ if $a(E, \cF,B,\bM)>0$ (resp. $\geq 0$) for any prime divisor $E$ exceptional over $X$ whose center on $X$ is the Zariski closure $\overline{\{P\}}$ of $P$, and we say $(X,\cF,B,\bM)$ is lc at $P$ if $a(E, \cF,B,\bM)\geq -\epsilon(E)$ for any prime divisor $E$ over $X$ whose center on $X$ is $\overline{\{P\}}$.

\subsection{Cone of curves} Let $K\subset \R^n$ be a strongly convex cone. We say a ray $R$ is \emph{exposed}  if there exists a linear functional $f\in (\R^n)^*$ such that $f(K)\geq 0$ and $R = K \cap \ker f$. 

Let $X$ be a normal projective variety and let $(X,\cF,B,\bM)$ be a gfq on $X$. Let $\NE(X)_{-\infty} $ be the subcone of $\NE(X)$ spanned by the images of $\NE(W)\to \NE(X)$ where $W$ is a non-lc center of $(X,\cF,B,\bM)$. Let $F$ be an extremal face of $\NE(X)$. A \emph{supporting function} of $F$ is a nef $\R$-divisor $H_F$ such that $F = \NE(X) \cap H_F^\perp $. We say $F$ is \emph{relatively ample at infinity} with respect to $(X,\cF,B,\bM)$ if 
$F \cap \NE(X)_{-\infty} = \{0\}$.

\subsection{Singularities in the sense of McQuillan} 

We recall the notion of singular locus for rank one foliations used by McQuillan \cite{McQ24}, following the notation of \cite[Section~4]{CS25a}. 

Let $X$ be a normal variety and let $\cF$ be a foliation of rank one on $X$ such that $K_\cF$ is $\Q$-Cartier. Let $x\in X$ be a point. We may take a local index one cover $\sigma: U' \to U$ associated to $K_\cF$ such that $U$ is an open neighborhood of $x$ and $\sigma^{-1}\cF$ is generated by a vector field $\partial$. We say $\cF$ is \emph{singular in the sense of McQuillan} at $x\in X$ if there exists an embedding $U' \to M$ to a smooth variety $M$ and a lift $\tilde{\partial}$ of $\partial$ on $M$ such that $\tilde{\partial}$ vanishes along $\sigma^{-1}(x)$. We denote by $\Sing^+\cF$ the locus of points $x\in X$ such that $\cF$ is singular in the sense of McQuillan at $x\in X$. 

\section{Birational modifications}

We begin with a discrepancy estimate for sequences of blowings-up along invariant centers over a foliated variety $(X,\cF)$ of rank one. The proof uses the same local calculation as in \cite[Lemma 1.1.3]{BM16}. 

\begin{lemma}\label{lemma:blowupdiscrep}
    Let $X$ be a normal variety and let $\cF$ be a foliation on $X$ of rank one. Assume that $K_\cF$ is Cartier. Let 
    \[
    h: Y = X_n \to X_{n-1}\to \cdots \to  X_0 = X
    \]
    be a sequence of blowing-ups $h_i : X_i \to X_{i-1}$ along $\cF$-invariant subvarieties over $X$. Assume that $Y$ is normal. Then for any $h$-exceptional prime divisor $E$ on $Y$, $a(E,\cF) \leq -\epsilon(E)$.
\end{lemma}

\begin{proof}
This is a local statement, and after shrinking $X$ to an affine open subset, we may assume that $X = \Spec A$ and $\cF$ is generated by $\partial \in \Der(A)$. Regard $\partial$ as a derivation of the function field $K(X)$. Let $Z_{i-1}\subset X_{i-1}$ be the blowing-up center of $h_i$.

First we show that $\partial$ lifts to a holomorphic vector field on $X_i$ for each $1\leq i \leq n$. Suppose that $\partial$ lifts to a holomorphic vector field on $X_{i-1}$. Since $Z_{i-1}$ is an $\cF$-invariant subvariety over $X$ and $\partial$ is holomorphic on $X_{i-1}$, we have $ \partial(\mathcal{I}) \subset \mathcal{I}$, where $\mathcal{I}$ is the ideal sheaf of $Z_{i-1}$.
Take an affine open subset $U = \Spec R$ of $X_{i-1}$, and let $I\subset R$ be the ideal associated to $Z_{i-1} \cap U$. Then $h_i^{-1}(U)$ is covered by $D_+(f) = \Spec R[I/f]$ for $f\in I \setminus \{0\}$, where 
\(
R[I/f]:= \{g/f^n\mid n\geq 0, g\in I^n\}.
\) 
Note that $\partial(I^n)\subset I^n$. For any $g \in I^n$, we have
\[
\partial \left(\frac{g}{f^n} \right) = \frac{f\partial g - ng\partial f}{f^{n+1}}\in R[I/f].
\]
Thus $\partial$ lifts to a holomorphic vector field on $X_i$.

Let $E$ be an $h$-exceptional prime divisor on $Y$ and let $t$ be a uniformizer of the DVR $\cO_{Y,\eta_E}$. Choose the smallest $j$ such that the center of $E$ on $X_j$ has codimension 1. By abuse of notation, we choose an open affine neighbourhood $U = \Spec R$
of the center of $E$ on $X_{j-1}$, and let $I\subset R$ be the ideal associated to $Z_{j-1}\cap U$, and choose a chart $\Spec R[I/f]\subset X_j$ for $f\in I \setminus \{0\}$, containing the center of $E$ on $X_j$ over $U$. We may write $f= t^d u$, where $d = \ord_E(f)$ and $u$ is a unit in $\cO_{Y,\eta_E}$. By the choice of $j$, we have $d>0$. Since $\partial f/f \in R[I/f] \subset \cO_{\eta_E}$, we have
\[
d\frac{\partial t}{t}
=
\frac{\partial f}{f}-\frac{\partial u}{u}
\in\cO_{\eta_E},
\]
which implies that $\ord_E(\partial t) \geq 1$.

Assume that $\lambda \in \Der(\cO_{\eta_E})$ generates the foliation $\cF_Y = h^{-1}\cF$ at $\eta_E$. We may write $\partial = t^{-a}v \lambda$ at $\eta_E$, where $v$ is a unit in $\cO_{\eta_E}$ and $a = a(E,\cF)$. If $a>0$, then $t^a \partial = v \lambda$. Since $\Der(\cO_{\eta_E})/\langle \lambda \rangle = (T_Y/T_{\cF_Y})_{\eta_E}$ is torsion-free, $\partial \in \langle \lambda \rangle$, and then $a\leq 0$, which leads to a contradiction. Thus $a\leq 0$. Moreover, if $E$ is not $\cF_Y$-invariant, then $\lambda(t)$ is a unit in $\cO_{\eta_E}$. Then 
    \[
    1\leq \ord_E(\partial t) = \ord_E(t^{-a}v\lambda(t)) = -a,
    \]
    and $a(E,\cF)= a \leq -1 = -\epsilon(E)$.
\end{proof}

Now we construct a log smooth birational model for an lc Gorenstein foliation of rank one such that the induced foliation remains lc, using the resolution algorithm given by Belotto \cite{Bel16}. One may regard it as a weak analogue of a foliated log resolution, comparing \cite[Lemma 6.2.2]{CHLX23} for algebraically integrable foliations or \cite[Lemma 3.3]{CS21} for co-rank one foliations on threefolds.

\begin{proposition}\label{prop:lcmodel}
    Let $X$ be a normal variety and let $\cF$ be a foliation of rank 1 on $X$. Assume that $K_\cF$ is Cartier, and $(X,\cF)$ is lc. Let $B$ be an effective $\R$-divisor and let $\bM$ be a b-$\R$-Cartier b-divisor on $X$. Then there exist a projective birational morphism $h:Y \to X$ and a reduced divisor $D$ on $Y$ with the following properties:
    \begin{enumerate}
        \item $(Y,D)$ is log smooth.
        \item $\Supp h^{-1}_* B \cup \Exc h \subset \Supp D$.
        \item $\bM$ descends to $Y$. 
        \item If $\cF_Y= h^{-1}\cF$, then $(Y,\cF_Y,D^{\ninv})$ is lc. 
    \end{enumerate}
\end{proposition}

\begin{proof}
    Since $K_\cF$ is Cartier, by the resolution algorithm given in \cite{BM97} or \cite[Theorem 1.1]{BM08}, there exists a projective birational morphism $f: X' \to X$, which is a sequence of blowing-ups along $\cF$-invariant centers over $X$, such that $X'$ is smooth, and the reduced exceptional divisor $E$ of $f$ has only simple normal crossings. Let $\cF' = f^{-1}\cF$. Since $(X,\cF)$ is lc, by Lemma~\ref{lemma:blowupdiscrep},  we have $K_{\cF'}+ E^{\ninv} = f^*K_\cF$, and $(X',\cF',E^{\ninv})$ is lc. 

    Take a projective birational morphism $\phi: W\to X'$ such that $\bM$ descends to $W$. By \cite[Theorem II.7.17]{Har77}, we may write $W= \Bl_\mathfrak{a} X'$ for a coherent ideal sheaf $\mathfrak{a}$ on $X'$. Let $B'=  \Supp f^{-1}_*B$, and $\mathcal{I}= \mathfrak{a}\cO_{X'}(-B')$. Using the notation in \cite{Bel16}, let $\theta= T_{\cF'}(-E^{\ninv})$ be a distribution on $X'$. Applying \cite[Theorem 1.1]{Bel16} to $(X',\theta, \mathcal{I},E)$, we obtain a projective birational morphism $g: Y \to X'$ and a reduced divisor $D$ on $Y$, such that the following hold:
    \begin{enumerate}
        \item[(i)] $(Y,D)$ is log smooth.
        \item[(ii)] $\Supp g^{-1}_* E \cup \Exc g \subset \Supp D$.
        \item[(iii)] $\mathcal{I} \cO_Y = \cO_Y(-\Delta)$, where $\Delta$ is an effective divisor such that $\Supp \Delta \subset \Supp D$.
        \item[(iv)] $(Y,\cF_Y,D^{\ninv})$ is lc, where $\cF_Y = g^{-1}\cF'$.
    \end{enumerate}
    Let $h= f\circ g: Y\to X$, then (1) and (4) follow. Note that $\mathcal{I} \cO_Y = (\mathfrak{a} \cO_Y) \otimes \cO_{Y}(-g^*B')$, we have $g^*B' \leq \Delta$, hence 
    \[
    \Supp h^{-1}_*B = \Supp g^{-1}_*B'\subset \Supp \Delta \subset \Supp D,
    \]
    which proves (2) together with (ii). Since $\mathfrak{a} \cO_Y$ is invertible, then $Y \to X$ factors through $W$. In particular, $\bM$ descends to $Y$, and (3) follows.
\end{proof}

\begin{corollary}\label{coro:Nlcclosed}
    Let $(X,\cF,B,\bM)$ be a gfq of rank one such that $K_\cF$ is $\Q$-Cartier. Then the non-lc locus $\Nlc(X,\cF,B,\bM)$ is Zariski closed in $X$.
\end{corollary}

\begin{proof}

    By \cite[Definition/Summary~I.ii.6]{MP13}, $\Nlc (X,\cF)$ is Zariski closed. Since $\Nlc(X,\cF)\subset \Nlc(X,\cF,B,\bM)$, we may replace $X$ with $X \setminus \Nlc(X,\cF)$ and assume that $(X,\cF)$ is lc.

    Let $x$ be a closed point of $X$ not contained in $\Nlc(X,\cF,B,\bM)$. It suffices to show that $(X,\cF,B,\bM)$ is lc over a Zariski open neighborhood of $x$. It is a local statement at $x$, so by \cite[Lemma 2.12]{Li25}, we may replace $X$ by an index one cover of $K_\cF$ and assume that $K_\cF$ is Cartier. Then let $h: Y\to X$ be a projective birational morphism and $D$ be a reduced divisor on $Y$ given by Proposition~\ref{prop:lcmodel}. We may write
    \[
    K_{\cF_Y} + B_Y + \bM_Y = h^*(K_\cF + B + \bM_X),
    \]
    where $\cF_Y = h^{-1}\cF$. Note that $\Supp B_Y \subset \Supp D$. Let $Z_Y = \Supp (B_Y - D^{\ninv})^+$. For any prime divisor $E$ contained in $\Supp Z_Y$, $\mu_E B_Y >\mu_E D^{\ninv} = \epsilon (E)$. Let $Z=h(Z_Y)$, then $Z\subset \Nlc(X,\cF,B,\bM)$. Furthermore, $B_Y \leq D^{\ninv}$ over $U = X\setminus Z$, then since $(Y,\cF_Y,D^{\ninv}) $ is lc, $(Y,\cF_Y,B_Y)$ is sub-lc over $U$. Since $\bM$ descends to $Y$, $(Y,\cF_Y,B_Y,\bM)$ is sub-lc over $U$, then $(X,\cF,B,\bM)$ is lc over $U$. Thus $\Nlc(X,\cF,B,\bM)=Z$ is Zariski closed in $X$.
\end{proof}

Next we consider rank one foliations that are terminal at every closed point. For $\Q$-Gorenstein rank one foliations, \cite[Lemma 2.9]{CS25b} gives a local characterization of this property. The following lemma shows that this property is preserved under finite group quotients.

\begin{lemma}\label{lemma:grpquotient}
    Let $X$ be a normal analytic space and let $\cF$ be a foliation of rank $1$ on $X$. Assume that $K_\cF$ is Cartier and $\cF$ is terminal at every closed point. Assume that there is a finite group $G$ acting on $X$ such that $\cF$ is $G$-invariant. Let $\cF_G$ be the induced foliation on $X/G$. Then $K_{\cF_G}$ is $\Q$-Cartier and $\cF_G$ on $X/G$ is terminal at every closed point.
\end{lemma}

\begin{proof}
    Let $y \in X/G$ be any closed point. We claim that there exists an analytic neighborhood $U$ of $y$ and a quasi-\'etale cover $q: V\to U $ such that $q^{-1}(\cF_G|_U)$ is induced by a submersion $f: V\to B$. 

    Let $\pi:X \to X/G$ be the quotient map and let $x\in \pi^{-1}(y)$ be a closed point. After replacing $G$ by the stabilizer of $x$, and then replacing by the quotient of $G$ by the kernel of the action, we may assume that $G$ acts faithfully on $X$ and fixes $x$. 
    By \cite[Lemma 2.9]{CS25b}, $x$ is not contained in $\Sing^+\cF$. Then by \cite[Lemma 1.2.1]{BM16}, there exists an analytic neighborhood $\tilde W$ of $x$ such that
    $\tilde W\cong\Delta_t\times Z$, where $Z$ is a normal irreducible analytic space, $x$ corresponds to $(0,z_0)$, and $\cF$ is generated by $\partial_t$.
    Let $p_Z: \tilde W\to Z$ be the projection.
    
    Since $\cF$ is $G$-invariant and $G$ fixes $x$, for any $g \in G$ there exists $\chi(g)\in \C^*$ such that $dg_x(\partial_t|_x) = \chi(g)\partial_t|_x$. Since $d(gh)_x = dg_x \circ dh_x$ for any $g,h \in G$, $\chi: G \to \C^*$ is a character. Let $W\subset \Tilde{W}$ be a $G$-invariant open neighborhood of $x$. Define 
    \[
    s :=\frac{1}{\lvert G \rvert} \sum_{g \in G} \chi(g)^{-1}g^*t \in \fm_{W,x}.
    \]
    Then we have $\frac{\partial s}{\partial t}(x)= (ds)_x(\partial_t|_x) = 1$ and $g^*s = \chi(g)s$ for any $g\in G$. By the holomorphic inverse function theorem, the map $(s,p_Z): W\to \C\times Z$ is a local analytic isomorphism at $x$. Since $\partial_t = \partial_t(s)\partial_s$ and $\partial_t(s)$ is a unit near $x$, $\cF$ is generated by $\partial_s$ near $x$. 
    
    Moreover, since $\cF$ is $G$-invariant,
    $\partial_s(g^*p_Z^*a)=0$ for every
    $a\in\cO_{Z,z_0}$ and every $g\in G$. 
    Thus for any $g\in G$, $g$ induces a local automorphism $g_Z$ of the germ $(Z,z_0)$ such that $p_Z \circ g = g_Z \circ p_Z$. Thus after shrinking $W$ to a $G$-invariant neighborhood of $x$, we may assume that $W \cong \Delta_s \times Z$ and $g(s,z)= (\chi(g)s, g_Z(z))$ for any $g\in G$ and $(s,z)\in \Delta_s \times Z$.

    Let $H$ be the subgroup of elements of $G$ acting trivially on $Z$ and $K$ be the subgroup of elements of $G$ acting trivially on $\Delta_s$. Then $H,K$ are normal subgroups of $G$, and $H\cap K = \{e\}$ as $G$ acts on $W$ faithfully. Since $H$ acts faithfully on $\Delta_s$, $\chi|_H$ is injective. Hence $H$ is cyclic. 
    Let $L= HK= H\times K$, $V = W/L$ and $r = \lvert H\rvert$. Then we have
    \[
    V = W/L \cong (\Delta_s/H)\times (Z/K) \cong \Delta_u \times Z/K,
    \]
    where $u = s^r$. Thus the induced foliation $\cF_L$ on $V$ is defined by the projection $V \to Z/K$. Moreover, $G/L$ acts faithfully on both $\Delta_u$ and $Z/K$. Thus for any $\gamma \in G/L$, if $\gamma$ is not the identity, the fixed locus of $\gamma$ on $V$ has codimension at least 2 in $V$. Let $U = W/G$. Then $q: V\to U$ is quasi-\'etale and $\cF_L = q^{-1}(\cF_G|_U)$. Let $m = \lvert G \rvert$. Then $(du)^{\otimes m}$ gives a nowhere vanishing $G/L$-invariant section of $\cO_V(mK_{\cF_L})$. Since $q$ is quasi-\'etale and $U$ is normal, $(du)^{\otimes m}$ induces a nowhere vanishing section of $\cO_U(mK_{\cF_G})$. In particular, $K_{\cF_G}$ is $\Q$-Cartier. Thus by \cite[Lemma 2.9]{CS25b}, $\cF_G$ is terminal at $y$. 
\end{proof}

Now we show that for a foliated variety $(X,\cF)$ of rank one such that $K_\cF$ is $\Q$-Cartier, if $\cF$ is terminal at every closed point, we may extract an $\cF$-invariant subvariety $W$ to an invariant prime divisor $E$ on $Y$ over $X$ such that the induced foliation $\cF_Y$ and the restricted foliation $\cF_{E^\nu}$ are both terminal at every closed point. The following proposition is stated in the analytic setting for later use.

\begin{proposition}\label{prop:extraction}
Let $X$ be a normal analytic space and let $\cF$ be a foliation of rank one on $X$. Assume that $K_\cF$ is $\Q$-Cartier and that $\cF$ is terminal at every closed point of $X$. Let $S$ be an $\cF$-invariant irreducible closed analytic subspace of $X$ with $S\neq X$. Then for any relatively compact open subset $U$ of $X$, there exists a projective birational morphism $h: Y \to U$ with the following properties:

\begin{enumerate}
    \item $Y$ has quotient singularities. In particular, $Y$ is $\Q$-factorial klt.
    \item Let $\cF_U = \cF|_U$ and $\cF_Y = h^{-1}(\cF_U)$. Then $\cF_Y$ is terminal at every closed point of $Y$, and $K_{\cF_Y} = h^*K_{\cF_U}$.
    \item For each irreducible component $S_0$ of $S\cap U$, there exists an $\cF_Y$-invariant prime divisor $E$ on $Y$ such that $h(E)= S_0$.
    \item For each $E$ as in (3), let $E^\nu\to E$ be the normalization map. Then $E^\nu$ has quotient singularities, and the restricted foliation $\cF_{E^\nu}$ of $\cF_Y$ on $E^\nu$ is terminal at every closed point of $E^\nu$.
\end{enumerate}
\end{proposition}

\begin{proof}
    Since $U$ is relatively compact in $X$, we may choose a finite cover by analytic open sets $U = \bigcup_{i=1}^m U_i$ such that $\cO_{U_i}(r_iK_{\cF}) \cong \cO_{U_i}$. Take the index-one covering Deligne-Mumford stack $\pi: \mathscr{U}\to U$ for $K_\cF$. Then $\mathscr{U}|_{U_i}\cong [V_i/\mu_{r_i}]$, and $\sigma_i: V_i \to U_i$ is the local cyclic index-one cover of $K_\cF$. Let $\cG_i = \sigma_i^{-1}\cF$, then $K_{\cG_i}$ is Cartier. By \cite[Lemma 2.8]{CS25b}, $\cG_i$ is terminal at every closed point of $V_i$. 
    
    By \cite[Lemma 2.9]{CS25b} and \cite[Lemma 1.2.1]{BM16}, we may assume that there is an open covering $V_i = \bigcup_j W_{ij}$ such that $W_{ij} \cong \Delta \times Z_{ij}$ where $\Delta $ is the unit disc, $Z_{ij}$ is a normal analytic space and $\cG_i|_{W_{ij}}$ is induced by the projection $p_{ij}:W_{ij} \to Z_{ij}$.

    Let $\mathscr{S} = (\pi^{-1}(S\cap U))_{\text{red}}$ and let $\mathscr{I}$ be the ideal of $\mathscr{S}$. For each $i$, let $S_i = (\sigma_i^{-1}(S\cap U_i))_{\text{red}}$. Since $S$ is $\cF$-invariant, there is a dense open subset $S^\circ\subset S$ such that $K_{\cF}$ is Cartier over $S^\circ$. Since $\sigma_i$ is \'etale over $S^\circ$, each component of $\sigma_i^{-1} (S^\circ\cap U_i)$ is $\cG_i$-invariant, hence each component of $S_i = \overline{\sigma_i^{-1} (S^\circ\cap U_i)}$ is also $\cG_i$-invariant. Then there exists a reduced analytic subspace $T_{ij} \subset Z_{ij}$ such that $S_i\cap W_{ij} \cong \Delta \times T_{ij}$.

    By \cite[Theorems 1.1 and 1.3]{BM08}, we may apply functorial principalization of $(\mathscr{U},\mathscr{I})$, and obtain $\tilde{h}: \mathscr{Y}\to \mathscr{U}$ such that $\mathscr{Y}$ is smooth and $\mathscr{I}\cO_{\mathscr{Y}} = \cO_{\mathscr{Y}}(-\mathscr{D})$, where $\mathscr{D}$ is a normal crossing divisor on $\mathscr{Y}$. Here the construction on each $V_i$ descends to $\mathscr{U}$ by the \'etale functoriality of the algorithm. Since the resolution sequence has finite length locally and $U$ is relatively compact, $\tilde{h}$ is a composition of finitely many blowing-ups. Let $q_i: \tilde{V}_i \to V_i$ be the base change of $\tilde{h}$ along $V_i \to \mathscr{U}$ and let $q_{ij}:\tilde{W}_{ij}\to W_{ij}$  be the base change of $\tilde{h}$ along $ W_{ij}\to V_i \to \mathscr{U}$. Then by the smooth functoriality of the algorithm, 
    \[
    q_{ij}= \id_\Delta \times r_{ij}: \tilde{W}_{ij} = \Delta \times \tilde{Z}_{ij} \to W_{ij} = \Delta \times Z_{ij},
    \]
    where $r_{ij}: \tilde{Z}_{ij} \to Z_{ij}$ is the principalization of $(Z_{ij}, I_{T_{ij}})$. Then $\tilde{\cG_i} = q_i^{-1}\cG_i$ is locally induced by the projection $\tilde{W}_{ij} \to \tilde{Z}_{ij}$, and every prime divisor on $\tilde{V_i}$ which is  $q_i$-exceptional or contained in the pullback of $\mathscr{D}$ is $\tilde{\cG_i}$-invariant. We also have $K_{\tilde{\cG_i}} = q_i^*K_{\cG_i}$.

    Let $\rho: \mathscr{Y}\to Y$ be the coarse moduli map and let $h: Y \to U$ be the induced morphism. We write $\tilde{h}$ as the composition
    \[
    \tilde{h}: \mathscr{Y} =  \mathscr{U}^N \to \cdots \to \mathscr{U}^1 \to \mathscr{U}^0 = \mathscr{U}
    \]
    of the blowing-ups $\tilde{h}_j: \mathscr{U}^j \to \mathscr{U}^{j-1}$, and let $\pi_j: \mathscr{U}^j\to U^j$ be the coarse moduli maps. Let $\mathscr{L}_j$ be the tautological $\tilde{h}_j$-ample line bundle on $\mathscr{U}^j$. Let $l = \lcm(r_1,r_2,\dots,r_m)$, then $\mathscr{L}_j^{\otimes l}$ has trivial stabilizer action on closed points. Therefore by \cite[Theorem 10.3]{Alp13}, there exists a line bundle $L_j$ on $U^j$ such that $\mathscr{L}_j^{\otimes l} = \pi_j^*L_j $. Since relative ampleness is preserved under finite surjective morphisms, $L_j$ is ample over $U^{j-1}$. Thus the morphisms $U^j \to U^{j-1}$ are projective, and so is $h$. Since $\mathscr{Y}$ is smooth, $Y$ has quotient singularities, and (1) follows. 

    Denote by $\tau_i: \tilde{V_i}\to Y$ the local quotient charts of $Y$. Let $\Gamma$ be a prime divisor on $\tilde{V_i}$ ramified over $Y$. Since $\sigma_i$ is quasi-\'etale, $\Gamma$ is exceptional over $V_i$ and hence $\tilde{\cG_i}$-invariant. Then by the foliated Riemann-Hurwitz formula \cite[Lemma 3.4]{Dru21},
    \[
    \tau_i^*K_{\cF_Y} = K_{\tilde{\cG_i}}- \sum_\Gamma \epsilon(\tau_i(\Gamma)) (r_\Gamma-1) \Gamma = K_{\tilde{\cG_i}} = q_i^*K_{\cG_i} =  q_i^*\sigma_i^*K_{\cF_U} = \tau_i^*h^*K_{\cF_U},
    \]
    hence $K_{\cF_Y} = h^*K_{\cF_U}$. By Lemma~\ref{lemma:grpquotient}, $\cF_Y$ is terminal at every closed point of $Y$, and (2) follows.

Let $\mathscr{E}$ be a prime divisor on $\mathscr{Y}$ contained in $\Supp \mathscr{D}$ such that the image of $\mathscr{E}$ on $U$ is $S_0$, and let $E = \rho(\mathscr{E})$, then $E$ is $\cF_Y$-invariant.
Let $\nu: E^\nu \to E$ be the normalization, and let $x$ be any closed point of $E^\nu$. Let $v\in \tilde{V}_i$ be a closed point over $\nu(x)$ and let $D$ be a local irreducible analytic branch
of $\tau_i^{-1}(E)$ at $v$ corresponding to $x$.
Since $\operatorname{Supp}\mathscr D$ has normal crossings, $D$ is smooth. Moreover, $D$ is $\widetilde{\cG}_i$-invariant,
and the restricted foliation $\cG_D$ on $D$ is regular. Let
$H_D$ be the stabilizer of the local branch $(D,v)$.
Then  $(E^\nu,x)\cong(D/H_D,\overline v)$ as analytic germs, where $\overline v$ is the image of $v$ on $D/H_D$. Hence $E^\nu$ has quotient singularities. Since $\cG_D$ is regular and $H_D$-invariant, by Lemma~\ref{lemma:grpquotient} the restricted foliation $\cF_{E^\nu}$ is terminal at $x$. Thus (3) and (4) follow.
\end{proof}

\section{Proof of the main theorems}

The following lemma extends \cite[Proposition 4.1]{Li25} from invariant curves to invariant subvarieties.

\begin{lemma}\label{lemma:invlccenter}
Let $(X,\cF,B,\bM)$ be a gfq of rank 1, and let $W$ be an irreducible $\cF$-invariant subvariety of $X$. Assume that $K_\cF$ is $\Q$-Cartier and $W \neq X$. If $(X,\cF,B,\bM)$ is lc at the generic point of $W$, then $W$ is an lc center of $(X,\cF)$. In particular, $W \not\subset \Supp B$ and $\bM$ descends to $X$ near the generic point of $W$.
\end{lemma}

\begin{proof}
    First we construct an lc place $E$ of $(X,\cF)$ centered on $W$. If $W$ is a prime divisor on $X$, then $0 \leq \mu_W B \leq \epsilon(W) = 0$ and we may conclude. Now we assume that $\codim_X W \geq 2$. Since $W$ is $\cF$-invariant, $K_\cF$ is Cartier in a neighborhood of the generic point of $W$. Let $\phi: Y \to X$ be the normalized blowing-up of $X$ along $W$ and let $\cF_Y = \phi^{-1}\cF$. Take $E$ to be a $\phi$-exceptional prime divisor on $Y$ which dominates $W$. By the proof of Lemma~\ref{lemma:blowupdiscrep} and \cite[\S 3, Theorem]{Sei66}, $a(E,\cF)\leq -\epsilon (E)$. Since $(X,\cF,B,\bM)$ is lc at the generic point of $W$, we have 
    \[
    -\epsilon (E) \leq a(E,\cF,B,\bM) \leq a(E,\cF) \leq -\epsilon(E),
    \]
    which implies that $W$ is an lc center for both $(X,\cF)$ and $(X,\cF,B,\bM)$. In particular, $W \not\subset \Supp B$. Let $f: X' \to X $ be a birational morphism such that $\bM$ descends to $X'$. We may assume $E$ is a prime divisor on $X'$. Let $P= f^*(B + \bM_X)-\bM_{X'}$. Then $-P$ is nef over $X$ and $f_*P = B \geq 0$, and by the negativity lemma, $P \geq 0$. Since $E$ is an lc place of both $(X,\cF)$ and $(X,\cF,B,\bM)$, $E \not\subset \Supp P$. By \cite[Lemma 3.39(2)]{KM98}, we have $W \not\subset f(\Supp P)$. Then there exists a neighborhood $U \subset X$ of the generic point of $W$ such that $B|_U=0$ and $\bM_{X'}|_{f^{-1}(U)} = (f|_U)^*(\bM_X|_U)$, which means $\bM|_U$ descends to $U$.  
\end{proof}

Next we show that a non-invariant compact curve is contained in an invariant surface if the foliation is terminal at every closed point of the curve.

\begin{lemma}\label{lemma:invsurf}
    Let $X$ be a normal analytic space of dimension $\geq 2$ and let $\cF$ be a foliation of rank $1$ on $X$ such that $K_\cF$ is $\Q$-Cartier. Let $C$ be an irreducible non-$\cF$-invariant compact curve such that $\cF$ is terminal at every closed point of $C$. Then after replacing $X$ with an analytic neighborhood of $C$, there is an irreducible $\cF$-invariant surface $S$ on $X$ which contains $C$.
\end{lemma}

\begin{proof}
Let $P\in C$ be any closed point. By \cite[Lemma 2.9]{CS25b}, there exists an analytic neighborhood $U$ of $P$ and a quasi-\'etale cover $\sigma: V\to U$ with Galois group $H$. Let $\cG = \sigma^{-1}\cF$, then $K_\cG$ is Cartier and $\cG$ is terminal at every closed point $Q\in \sigma^{-1}(C\cap U)$. By \cite[Lemma 1.2.1]{BM16}, we may take an analytic neighborhood $W_Q$ of $Q$ such that $W_Q \cong \Delta \times Z_Q$, where $\Delta $ is the unit disc and $Z_Q$ is normal, such that $\cG$ is induced by the projection $W_Q \to Z_Q$ near $Q$. Since $C$ is not $\cF$-invariant, there exists a $\cG$-invariant surface germ $\Gamma_Q\subset W_Q$ such that $\Gamma_Q$ is the minimal invariant germ containing $\sigma^{-1}(C)$ around $Q$. In particular, for every $g \in H$ we have $\Gamma_{g(Q)}=g(\Gamma_Q)$ as surface germs. Define $S_P = \sigma(\Gamma_Q)$ where $Q\in \sigma^{-1}(P)$ is a closed point. Thus $S_P$ is independent of the choice of $Q$, and $S_P$ is the minimal invariant surface germ containing $C$ around $P$. Since $C$ is compact, after replacing $X$ with an analytic neighborhood of $C$, the germs $S_P$ patch together and induce an irreducible $\cF$-invariant surface $S$ containing $C$.
\end{proof}

Next we consider the adjunction of singularities. 
\begin{lemma}\label{lemma:noninvcurveadjLC}
    Let $X$ be a $\Q$-factorial analytic space and let $(X,\cF,B,\bM)$ be a rank 1 gfq on $X$. Let $C$ be a non-$\cF$-invariant curve on $X$ not contained in $\Sing^+\cF$, and let $S$ be an $\cF$-invariant prime divisor such that $S$ is not contained in $\Supp B$ and $C\subset S$. Let $\nu:S^\nu\to S$ be the normalization map. Then we have
    \[
    (K_\cF + B + \bM_X)|_{S^\nu} = K_{\cF_S} + B_S + \bM^S_{S^\nu},
    \]
    where $\cF_S$ is the restricted foliation of $\cF$ on $S^\nu$, $B_S$ is an effective $\R$-divisor on $S^\nu$, and $\bM^S = \bM|_{S^\nu}$ is the restricted b-divisor. Moreover, if $(X,\cF,B,\bM)$ is lc at the generic point of $C$, then $(S^\nu,\cF_S,B_S,\bM^S)$ is lc at the generic point of $\nu^{-1}(C)$.
\end{lemma}

\begin{proof}

By \cite[Proposition 2.13]{Li25}, we have $B_S$ is an effective $\R$-divisor on $S^\nu$. Then the lemma follows from the same proof as in \cite[Lemma 4.5]{Li25} with the same notation. Here the argument in \cite[Lemma 4.5]{Li25} applies in arbitrary dimension.
\end{proof}

The following proposition characterizes irreducible components of the null locus of a supporting function of an exposed face, which is a key step of the proof of the cone theorem.

\begin{proposition}\label{prop:nulllocus}
    Let $X$ be a normal projective variety, and let $(X,\cF,B,\bM)$ be a gfq of rank 1. Assume that $K_\cF$ is $\Q$-Cartier.
    Let $A$ be an ample $\R$-divisor such that $H = K_\cF + B + \bM_X + A$ is nef. Then for any irreducible component $W$ of $\Null H$ such that $W$ is not contained in the non-lc locus $\Nlc(X,\cF,B,\bM)$, either $W$ is contained in $\Sing^+\cF$, or $W$ is $\cF$-invariant.
\end{proposition}

\begin{proof}
    Assume that $W$ is not contained in $\Sing^+\cF$ and $W$ is not $\cF$-invariant. Then $W \neq X$ and $H$ is big. Let $V \to W$ be the normalization map. Let $D$ be a sufficiently small ample $\R$-divisor such that $A - D$ is ample and $H - D$ is a big $\Q$-Cartier divisor. By \cite[Theorem 1.4]{birkar_augmented_2017}, we may assume that $\Null(H) = \mathbf{B}(H-D)$. Let $0 \leq G \sim_\Q H - D$ be an effective $\Q$-Cartier divisor. Since $H|_V$ is not big, $G|_V$ is not pseudo-effective. Then there exists a general curve $C_V \subset V$ such that $ G|_V\cdot C_V<0$. Let $C$ be the image of $C_V$ on $X$, then $G \cdot C<0$. Since $C_V$ is general on $V$, we may assume $C$ is not $\cF$-invariant and $C \not\subset \Sing^+\cF$. Since $W$ is not contained in the non-lc locus $\Nlc(X,\cF,B,\bM)$, by Corollary~\ref{coro:Nlcclosed}, we may assume $(X,\cF,B,\bM)$ is lc at the generic point of $C$. Set $K = K_\cF + B +\bM_X$. We have
\[
K\cdot C = (H -D)\cdot C -(A-D)\cdot C = G|_V \cdot C_V - (A-D)\cdot C <0.
\]

By \cite[Lemma 4.6]{CS25a}, after possibly replacing $X$ by an analytic neighborhood of $C$, there exists a birational morphism $p:X'\to X$ such that the following statements hold:

\begin{itemize}
    \item Let $\cF' = p^{-1}\cF$, and we may write
    \[
    K_{\cF'} + B' + \bM_{X'} = p^*K,
    \]
    where $B'$ is an effective $\R$-divisor on $X'$.
    \item $p^{-1}$ is an isomorphism at the generic point of $C$. Let $C'$ be the strict transform of $C$ on $X'$, then $\cF'$ is terminal at every closed point of $C'$.
    \item There exists an $\cF'$-invariant surface $\Gamma$ containing $C'$.
\end{itemize}

Moreover, by the construction in \cite[Lemma 4.6]{CS25a} and Lemma~\ref{lemma:grpquotient}, we have $K_{\cF'}$ is $\Q$-Cartier in a neighborhood of $C'$. Let $x$ be a general point of $C$ such that $x$ is contained in no other component of $\Null H$ but $W$. Since $p^{-1}$ is an isomorphism near $x$, and $W$ is not $\cF$-invariant, then $p(\Gamma) $ is not contained in $W$, and moreover, $\Gamma$ is not contained in $p^{-1}(\Null H)$. Thus we may take general $0 \leq G\sim_\Q H -D$ such that $\Gamma$ is not contained in $\Supp p^*G$. After replacing $X,\cF,B,C,G$ with $X',\cF',B',C',p^*G$ respectively, and possibly replacing $X$ with a neighborhood of $C$, we may assume the following properties hold:
\begin{enumerate}
    \item $X$ is a normal analytic space, $(X,\cF,B,\bM)$ is a gfq of rank one and $K_\cF$ is $\Q$-Cartier.
    \item $C$ is a non-$\cF$-invariant irreducible compact curve on $X$, and there exists an irreducible $\cF$-invariant surface $\Gamma$ containing $C$.
    \item $\cF$ is terminal at every closed point of $C$.
    \item $(X,\cF,B,\bM)$ is lc at the generic point of $C$.
    \item $K\cdot C = (K_\cF + B + \bM_X)\cdot C<0$.
    \item There exists an effective $\Q$-Cartier divisor $G \geq 0$ such that $G \cdot C<0$ and $\Gamma$ is not contained in $\Supp G$.
\end{enumerate}

We prove by induction on dimension that there is no $(X,\cF,B,\bM),C,\Gamma$ and $G$ satisfying (1)--(6).
First we assume that $\dim X \geq 3$. 
By \cite[Lemma 2.9]{CS25b}, after replacing $X$ with an analytic neighborhood of $C$, we may assume that $\cF$ is terminal at every closed point of $X$. By Proposition~\ref{prop:extraction}, after replacing $X$ with a relatively compact open neighborhood of $C$, there exists a projective birational morphism $h:Y \to X$ and an $\cF_Y$-invariant prime divisor $E$ on $Y$ such that $h(E)= \Gamma$, where $\cF_Y = h^{-1}\cF$. 
Then we may write 
\[
K_{\cF_Y}+B_Y + \bM_Y = h^*K.
\]
Since $K_{\cF_Y} = h^*K_\cF$ and by the negativity lemma, we have $B_Y\geq 0$. Since $(X,\cF,B,\bM)$ is lc at the generic point of $C$ and $E$ is an $\cF_Y$-invariant prime divisor whose center on $X$ contains $C$, a local computation as in \cite[Remark 2.3]{CS21} shows that 
$E$ is not contained in $\Supp B_Y$. Let $S  \to E$ be the normalization map, and denote by $h_S:S\to X$ the composition of $\iota:S \to Y$ and $h$. Let $\cF_S$ be the restricted foliation of $\cF_Y$ on $S$. By \cite[Proposition 2.13]{Li25}, there exists an $\R$-divisor $B_S \geq 0$ such that
\[
K_{\cF_S}+ B_S + \bM^S_S = (K_{\cF_Y} + B_Y + \bM_Y)|_S,
\]
where $\bM^S = \bM|_S$ is the restricted b-divisor.

Since $h$ is projective, $h_S$ is also projective. Let $F$ be an irreducible component of $h_S^{-1}(C)$ which dominates $C$. We claim that $F$ is not $\cF_S$-invariant. Otherwise $\iota(F)$ is $\cF_Y$-invariant. Since $\cF_Y$ is terminal at every closed point of $Y$, then by \cite[Lemma 2.9]{CS25b}, $\cF_Y$ is canonical. Thus by \cite[Lemma 1.1.3]{BM16}, there exists a prime divisor $\Lambda$ over $Y$ such that $a(\Lambda, \cF_Y) = 0$ and the center of $\Lambda$ on $Y$ is $\iota(F)$, hence $a(\Lambda,\cF) = 0$ and the center of $\Lambda$ on $X$ is $C$. This contradicts \cite[Lemma 4.3]{CS25a}, which implies that $\cF$ is terminal at the generic point of $C$. 

Since $h_S$ is projective, by cutting $F$ by general hyperplane sections, we may find a general curve $C_S \subset F$ which dominates $C$. By Proposition~\ref{prop:extraction}, $\cF_S$ is terminal at every closed point of $S$. Since $C$ is not $\cF$-invariant and $C_S$ is general in $F$, we may take $C_S$ to be non-$\cF_S$-invariant. Let $C_Y = \iota(C_S)$, then since $(X,\cF,B,\bM)$ is lc at the generic point of $C$, $(Y,\cF_Y,B_Y,\bM)$ is lc at the generic point of $C_Y$. Thus by Lemma~\ref{lemma:noninvcurveadjLC}, $(S,\cF_S,B_S,\bM^S)$ is lc at the generic point of $C_S$. Denote by $d$ the degree of $C_S \to C$. Then
\[
\left(K_{\cF_S} + B_S + \bM^S_S\right)\cdot C_S = h_S^*K \cdot C_S = d K \cdot C<0.
\]
By Lemma~\ref{lemma:invsurf}, after replacing $S$ with an analytic neighborhood of $C_S$, we may assume that there exists an $\cF_S$-invariant surface $\Gamma_S$ on $S$ such that $C_S \subset \Gamma_S$. Let $G_S=h_S^*G$. Since $\Gamma$ is not contained in $\Supp G$, $E$ is not contained in $\Supp h^*G$, and hence $G_S\geq 0$, and $G_S \cdot C_S = d G\cdot C <0$.

Now we claim that $\Gamma_S$ is not contained in $\Supp G_S$. Otherwise, near a general point of $C_S$, we have $C \subset h_S(\Gamma_S) \subset \Supp G \cap \Gamma$. Since $\Gamma$ is not contained in $\Supp G$, we have $\Supp G \cap \Gamma = C$ near a general point of $C$. Thus $\Gamma_S\subset F$, which contradicts the fact that $C_S$ is general in $F$ and $F$ is not $\cF_S$-invariant. 

Therefore, $(S,\cF_S,B_S,\bM^S),C_S,\Gamma_S,G_S$ also satisfy (1)--(6) above. Thus we may apply induction on dimension and reduce to the case when $\dim X = 2$ under the setting of (1)--(6). 

Now we assume that $\dim X = 2$. By (4), we have $\mu_C B \leq 1$. In the following we may consider  Mumford's intersection theory on normal surfaces; see \cite[Section 3]{Li25} for example. By (6), we have $C^2<0$. Let $t = 1- \mu_C B \geq 0$. By \cite[Proposition 3.7(1)]{Li25}, we have $(K_\cF + B + tC + \bM_X)\cdot C \geq 0$. Note that although \cite[Proposition 3.7(1)]{Li25} is stated for projective surfaces, the proof also works in our case. However,
\[
 (K_\cF + B +tC + \bM_X)\cdot C = (K_\cF + B + \bM_X)\cdot C + tC^2<0,
\]
which leads to contradiction.
\end{proof}

\begin{remark}\label{rmk:HirSur}
The ample divisor in Proposition~\ref{prop:nulllocus} cannot be omitted. If $(X,\cF,B)$ of rank one is lc and $H = K_\cF+ B $ is nef, $\Null H$ may contain a non-invariant component not contained in $\Sing^+\cF$, as the following example shows. Let $n \geq 1$. Take $X = \mathbb{F}_n = \mathbb{P}_{\mathbb{P}^1} (\cO \oplus \cO(n))$ be the Hirzebruch surface, and let $\cF$ be the foliation induced by $\mathbb{F}_n \to \mathbb{P}^1$. Let $h$ be the negative section on $X$ and $f$ be a general fiber of $X \to \mathbb{P}^1$, hence $h^2=-n$, $f^2 = 0$ and $h\cdot f = 1$. Note that $\lvert h + nf\rvert $ is base-point-free, we may take general members $l_i\in \lvert h + nf\rvert$, $1\leq i \leq 4$, such that any three of $l_i$ do not pass through a common point. Let $B = h+ \frac{1}{2}(l_1 + l_2 + l_3 +l_4)$. Then by \cite[Lemma 4.4]{CS25a}, $(X,\cF,B)$ is lc. However, $K_\cF + B \sim_\Q h + nf$, and $\Null (K_\cF + B) = h$ is not $\cF$-invariant. This example also provides a counterexample to \cite[Corollary 4.9]{CS25a} as stated.
\end{remark}

Now we are ready to prove Theorem~\ref{MainThm:ConeTHM}.

\begin{proof}[Proof of Theorem~\ref{MainThm:ConeTHM}]

Set $K= K_\cF+B+\bM_X$. Let $R$ be a $K$-negative exposed ray which is relatively ample at infinity with respect to $(X, \cF, B, \bM)$, and let $H_R$ be a supporting function to $R$. In particular, $H_R$ is nef. By \cite[Lemma 8.2.7]{CHLX23}, we may assume that there exists an ample $\R$-divisor $A$ such that $H_R = K + A$. 

First we show that there exists an irreducible component $W$ of $\Null H_R$ such that $R$ is contained in the image of $\NE(W)\to \NE(X)$. By \cite[Theorem 1.4]{birkar_augmented_2017}, we may choose a sufficiently small ample $\R$-divisor $D$ such that $A-D$ is ample, $H_R-D$ is $\Q$-Cartier and $\Null H_R = \bB (H_R-D)$. Since $\NE(X)$ is the closed convex hull of the classes of irreducible curves on $X$ and $R$ is exposed, there exists a sequence of irreducible curves $\{\xi_i\}_{i=1}^\infty$ such that $R$ is the limit of rays $\R_{\geq 0}[\xi_i]$. Since $(H_R-D)\cdot R<0$, we have $(H_R -D)\cdot \xi_i<0$ for all sufficiently large $i$. In particular, $\xi_i$ is contained in $\bB(H_R - D)$ for sufficiently large $i$. After replacing $\{\xi_i\}$ by a subsequence, we may assume that each $\xi_i$ is contained in a fixed irreducible component $W$ of $\bB(H_R-D)$. Thus $R$ is contained in the image of $\NE(W)\to \NE(X)$. 
Let $V\to W$ be the normalization map. Then $H_R|_V$ is not big. Since $R$ is relatively ample at infinity, $W$ is not contained in the non-lc locus $\Nlc(X,\cF,B,\bM)$. By Proposition~\ref{prop:nulllocus}, either $W$ is contained in $\Sing^+\cF$, or $W$ is $\cF$-invariant.

If $W$ is contained in $\Sing^+\cF$, then $\xi_i\subset \Sing^+\cF$ and  
\[
K\cdot\xi_i = (H_R -D)\cdot\xi_i - (A-D)\cdot \xi_i <0.
\]
By \cite[Corollary 4.3]{Li25}, each $\xi_i$ is contained in the non-lc locus $\Nlc(X,\cF,B,\bM)$. By Corollary~\ref{coro:Nlcclosed}, after replacing $\{\xi_i\}$ by a subsequence, we may assume each $\xi_i$ is contained in a fixed component $W'$ of $\Nlc(X,\cF,B,\bM)$. Thus $R$ is contained in the image of $\NE(W')\to \NE(X)$, which contradicts our choice of $R$. Thus we may assume that $W$ is not contained in $\Sing^+\cF$ and is $\cF$-invariant.

Let $d = \dim W$. Since $H_R|_V$ is not big, there exists an integer $0\leq \nu \leq d-1$ such that 
\[
(H_R|_V)^{\nu+1}\cdot (A|_V)^{d-\nu-1} = 0, \quad\text{and } (H_R|_V)^{\nu}\cdot (A|_V)^{d-\nu} > 0.
\]
Let $D_i := H_R|_V$ for $1\leq i \leq \nu+1$, $D_j:=A|_V$ for $\nu+1<j\leq d$. Thus
\begin{align*}
    -K|_V \cdot D_2\cdot  \cdots \cdot  D_d &= (A-H_R)|_V \cdot (H_R|_V)^{\nu}\cdot (A|_V)^{d-\nu-1} \\
    &= (H_R|_V)^{\nu}\cdot (A|_V)^{d-\nu}-(H_R|_V)^{\nu+1}\cdot (A|_V)^{d-\nu-1}>0.
\end{align*}

If $W = X$, then since $\bM_X$ is pseudo-effective, $\bM_X \cdot D_2 \cdot \cdots D_d \geq 0$. If $W \neq X$, then
by Lemma~\ref{lemma:invlccenter}, there exists a prime divisor $T$ over $X$, whose center on $X$ is $W$, such that $a(T,\cF) = a(T,\cF,B,\bM)= -\epsilon(T)$. Let $f: Y \to X$ be a birational morphism such that $Y$ is smooth, $\bM$ descends to $Y$ and $T$ is a divisor on $Y$. We may write $\bM_Y + P = f^*\bM_X$ for $P\geq 0$ exceptional over $X$. Then $T \not\subset \Supp P$. Take a general subvariety $\Sigma \subset T$ such that $\Sigma \to W$ is a surjective, generically finite morphism. Since $\Sigma$ is general in $T$, we may assume that $\Sigma$ is not contained in $\Supp P$. Let $\Sigma'\to \Sigma$ be the normalization map, and write $q:\Sigma'\to V$ for the induced morphism. Therefore,
\begin{align*}
    \bM_X|_V \cdot D_2\cdot  \cdots \cdot  D_d &= \frac{1}{r} q^*(\bM_X|_V) \cdot q^*D_2 \cdot \cdots \cdot q^*D_d \\
    &= \frac{1}{r} (f^*\bM_X)|_{\Sigma'}\cdot q^*D_2 \cdot \cdots \cdot q^*D_d \\
    &= \frac{1}{r} (\bM_Y + P)|_{\Sigma'}\cdot  q^*D_2 \cdot \cdots \cdot q^*D_d \geq 0,
\end{align*}
where $r$ is the degree of $q$, and the last inequality holds as $\bM_Y|_{\Sigma'}$, $q^*D_i$ are nef on $\Sigma'$ and $P|_{\Sigma'}$ is effective. Then we have $-(K_\cF + B)|_V \cdot D_2\cdot \cdots \cdot D_d >0$. By \cite[Proposition-Definition 3.12]{CS25a}, there exists an $\R$-divisor $B_V \geq 0$ such that 
\[
K_{\cF_V} + B_V = (K_\cF + B)|_V,
\]
where $\cF_V$ is the restricted foliation of $\cF$ on $V$. By \cite[Corollary 2.28]{Spi20}, through any general point of $V$, there exists a curve $\xi_V$ tangent to $\cF_V$ with $H_R|_V \cdot \xi_V = 0$. Let $\xi$ be the image of $\xi_V$ in $X$, then $[\xi]$ generates $R$. By the proof of \cite[Lemma 4.7]{Li25}, $\xi$ is a rational curve with $0<-K\cdot \xi\leq 2$.

Therefore, for any $K$-negative exposed ray $R$ which is relatively ample at infinity with respect to $(X,\cF,B,\bM)$, there exists a rational curve $C$ tangent to $\cF$ such that $0 <-K\cdot C\leq 2$ and $R$ is generated by $[C]$. Let $\{R_j\}_{j\in \Lambda}$ be the set of $K$-negative extremal rays which are relatively ample at infinity. Standard convex geometry arguments show that for any ample $\R$-divisor $A$, the set
\[
\Lambda_A :=\{j\in \Lambda \mid R_j \subset \NE(X)_{K+ A<0}\}
\]
is finite, and for each $j \in \Lambda$, $R_j$ is an exposed ray. In particular, $\Lambda$ is at most countable, and we may conclude.
\end{proof}

As a direct application of Theorem~\ref{MainThm:ConeTHM}, we prove Theorem~\ref{MainThm:MMPoutputs}, showing that the MMP outputs are unique up to isomorphism for generalized foliated quadruples of rank one.

\begin{proof}[Proof of Theorem~\ref{MainThm:MMPoutputs}]
    Since $Y_1,Y_2$ are outputs of MMPs, $Y_1$ and $Y_2$ are both normal projective $\Q$-factorial varieties. Let $g: W \to X$ be a birational morphism which resolves the indeterminacy of $f_1$ and $f_2$ and let $p: W \to Y_1$ and $q: W \to Y_2$ be the induced birational morphisms. Let $\phi = f_2\circ f_1^{-1}: Y_1 \bir Y_2$. We may assume that $\bM$ descends to $W$. Let $\cF_1$, $\cF_2$ and $\cG$ be the foliations on $Y_1$, $Y_2$ and $W$ induced by $\cF$ respectively. For each $i = 1,2$, let $B_i$ be the strict transform of $B$ on $Y_i$. Then $K_{\cF_i} + B_i + \bM_{Y_i}$ is nef. By the negativity lemma, there exists an $\R$-divisor $E\geq 0$ on $W$ such that 
    \[
    g^*(K_\cF + B + \bM_X) = p^*(K_{\cF_1}+ B_1 + \bM_{Y_1}) + E = q^*(K_{\cF_2}+ B_2 + \bM_{Y_2}) + E,
    \]
    and $E$ is exceptional over both $Y_1$ and $Y_2$. Since $f_1$ and $f_2$ are both $(K_\cF + B + \bM_X)$-MMPs, a prime divisor on $X$ is contracted by $f_i$ if and only if its strict transform on $W$ is contained in $\Supp E$. Thus $\phi$ is small.

    Let $A$ be a very ample divisor on $Y_2$, let $\bN =\overline{A}$ and let $H = \phi^{-1}_*A = \bN_{Y_1}$. By the negativity lemma, $F = p^*H -q^*A \geq 0$. We claim that $\Supp F \subset \Supp E$. Otherwise let $D$ be a prime divisor contained in $\Supp F$ but not in $\Supp E$. Then $D$ is exceptional over $Y_1$, and \[
    a(D,\cF,B,\bM) = a(D,\cF_1,B_1,\bM).
    \]
    By the same argument as the proof of \cite[Lemma 3.38]{KM98}, both $f_1$ and $f_2$ are isomorphic near the generic point of the center of $D$ on $X$. In particular, $p^*H = q^* A$ near the generic point of $D$, which contradicts $D \subset \Supp F$, and the claim holds. Thus we may take $t>0$ such that $E -tF\geq 0$.

    Let $B_W$ be the $\R$-divisor on $W$ satisfying
    \[
    K_\cG + B_W + \bM_W = g^*(K_\cF + B + \bM_X).
    \]
    Then $(W,\cG,B_W,\bM)$ is sub-lc. Since $\bN$ descends to $Y_2$, $\bN$ also descends to $W$, and $(W,\cG,B_W,\bM+t\bN)$ is sub-lc. Note that
    \[
    K_\cG + B_W -E + tF + \bM_W +t\bN_W = p^*(K_{\cF_1} + B _1+ \bM_{Y_1} + t\bN_{Y_1}).
    \]
    Since $E -tF\geq 0$, $(W,\cG,B_W-E+tF,\bM + t \bN)$ is sub-lc, then $(Y_1,\cF_1,B_1,\bM + t\bN)$ is lc. 

    We claim that $K_{\cF_1} + B_1 + \bM_{Y_1} + t\bN_{Y_1}$ is nef. Otherwise by Theorem~\ref{MainThm:ConeTHM}, there exists an $\cF_1$-invariant curve $C$ on $Y_1$ such that $(K_{\cF_1} + B_1 + \bM_{Y_1} + t\bN_{Y_1}) \cdot C <0$. Since $K_{\cF_1}+ B_1 + \bM_{Y_1}$ is nef, we have $\bN_{Y_1} \cdot C <0$. Since $(Y_1,\cF_1,B_1,\bM + t\bN)$ is lc, $(Y_1,\cF_1,B_1,t\bN)$ is also lc. Then \cite[Proposition 4.1]{Li25} implies that $t\bN_{Y_1}\cdot C \geq 0$, which leads to a contradiction. Note that
    \[
    tF = p^*(K_{\cF_1} + B_1 + \bM_{Y_1} + tH)-q^*(K_{\cF_2}+ B_2 + \bM_{Y_2} + tA)
    \]
    is $q$-exceptional and nef over $Y_2$. By the negativity lemma, $F \leq 0$, hence $F = 0$, or equivalently, $p^*H = q^*A$. By the rigidity lemma, $\phi$ is a morphism. Since $\phi$ is a projective small birational morphism between $\Q$-factorial varieties, $\phi$ must be an isomorphism.
\end{proof}

\begin{remark}\label{rmk:MMPoutputs}
    The example given by \cite[Remark 3.4]{CMS26} shows that rank one foliated minimal models are not unique in birational class. We may modify the example there, showing that the uniqueness also fails for $K_\cF$-minimal models of $(X,\cF)$. Let $E$ be an elliptic curve, and let $Y \bir Y'$ be a nontrivial flop of smooth projective threefolds. Define $X = E \times Y$ and $X' = E \times Y'$, and let $\cF$ (resp. $\cF'$) be the foliation induced by the projection $X \to Y$ (resp. $X' \to Y'$). Then $\phi: X \bir X'$ is small, and $K_\cF$ and $K_{\cF'}$ are both linearly trivial. Then $\phi$ and $\id_X$ are both $K_\cF$-negative, hence are $K_\cF$-minimal models of $(X,\cF)$, while $\phi:X\bir X'$ is not an isomorphism. 
\end{remark}

\bibliographystyle{alpha}
\bibliography{ref-cone}

\end{document}